\documentclass[11pt]{amsart}
\usepackage[a4paper,margin=28mm]{geometry}
\usepackage{amsmath,amssymb,amsfonts,amsthm,mathrsfs}
\usepackage{xcolor}
\usepackage[colorlinks=true,linkcolor=blue!55!black,citecolor=blue!55!black,urlcolor=blue!55!black]{hyperref}
\numberwithin{equation}{section}
\newtheorem{theorem}{Theorem}[section]
\newtheorem{maintheorem}{Theorem}

\newtheorem{proposition}[theorem]{Proposition}
\newtheorem{lemma}[theorem]{Lemma}
\newtheorem{corollary}[theorem]{Corollary}
\theoremstyle{definition}\newtheorem{definition}[theorem]{Definition}
\theoremstyle{remark}
\newcommand{\one}{\mathbf1}
\newcommand{\norm}[1]{\lVert#1\rVert}
\newcommand{\ip}[2]{\langle#1,#2\rangle_{L^2(\mu)}}
\newcommand{\cE}{\mathcal E}
\title[Nonexistence of fractional Lane--Emden supersolutions]{Volume criteria for nonexistence of positive fractional Lane--Emden supersolutions}
\author{Haoxuan Cheng}
\address{School of Mathematical Sciences, Fudan University, Shanghai, China}
\email{hxcheng25@m.fudan.edu.cn}
\date{}
\begin{document}
\raggedbottom
\begin{abstract}
For $0<s<1$ and $q>1$, we prove that every nonnegative weak supersolution of $(-\Delta)^su\ge u^q$ on a complete connected Riemannian manifold without boundary is zero if $\int_1^\infty R^{2sq-1}V(R)^{1-q}\,dR=\infty$, where $V(R)$ is geodesic ball volume. On an infinite connected locally finite loopless weighted graph, the corresponding volume series excludes nonzero nonnegative pointwise supersolutions for the fractional normalized Laplacian, with balls defined by the original combinatorial distance. Both conclusions also hold for finite Green integral supersolutions.
\end{abstract}
\maketitle

\section{Introduction}
Volume-growth criteria for semilinear elliptic inequalities use the volumes of metric balls to exclude nonzero nonnegative supersolutions. We study this question on complete Riemannian manifolds and infinite weighted graphs for the inequality
\[
 A^su\ge u^q,\qquad 0<s<1,\quad q>1.
\]
Here $A$ is the nonnegative Laplace--Beltrami operator or the normalized graph Laplacian, and $A^s$ is its spectral fractional power. In either setting, $V(R)$ denotes the volume of a ball of radius $R$ about a fixed center. The aim is to determine volume conditions that force every nonnegative supersolution to vanish, including in geometries where global Green-kernel estimates are unavailable.

For the ordinary Laplacian on manifolds, Grigor'yan and Kondratiev \cite{GrigoryanKondratiev2010} developed capacity and volume tests for semilinear inequalities. Grigor'yan and Sun \cite{GrigoryanSun2014} proved nonexistence under the critical bound
\[
 V(R)\le CR^{2q/(q-1)}(\log R)^{1/(q-1)}
\]
and established sharpness of its power and logarithmic exponents.

Green potentials give another approach to the same problem. Kalton and Verbitsky \cite{KaltonVerbitsky1999} related superlinear integral equations to weighted norm and capacity inequalities. Grigor'yan and Verbitsky obtained pointwise estimates for elliptic inequalities \cite{GrigoryanVerbitsky2019} and nonlinear integral inequalities \cite{GrigoryanVerbitsky2020}. Grigor'yan, Sun, and Verbitsky \cite{GrigoryanSunVerbitsky2020} derived manifold Green-function criteria that, under volume doubling and suitable two-sided Green estimates, reduce to convergence of a volume integral with exponent $2q-1$. Heat-kernel estimates supply the geometric information needed in this reduction. The work of Grigor'yan \cite{Grigoryan1992} and Saloff-Coste \cite{SaloffCoste1992} connects volume and Poincar\'e inequalities to Harnack inequalities and Gaussian bounds; Coulhon and Sikora \cite{CoulhonSikora2008} obtain Gaussian upper bounds from Davies--Gaffney estimates and on-diagonal bounds.

The extension approach of Caffarelli and Silvestre \cite{CaffarelliSilvestre2007}, and its spectral formulation by Stinga and Torrea \cite{StingaTorrea2010}, express a fractional operator through a local equation in an additional variable. Kwa\'snicki \cite{Kwasnicki2017} compares the Euclidean definitions of the fractional Laplacian. Gu, Huang, and Sun \cite{GuHuangSun2026} used fractional Green potentials to characterize the existence of positive Green integral supersolutions on manifolds satisfying volume doubling, a scale-invariant Poincar\'e inequality, and fractional transience. Their criterion is convergence of the volume integral appearing in the following theorem. We prove its nonexistence implication on every complete connected manifold without boundary, including weak supersolutions. The precise weak and Green integral solution classes used in both of our theorems are given in Section~\ref{sec:upper}.

\begin{samepage}
\begin{maintheorem}\label{thm:intro-manifold}\label{thm:manifold}
Let $M$ be a smooth complete connected Riemannian manifold without boundary, of positive dimension, and let $A=-\Delta$ be its nonnegative Friedrichs Laplacian. Fix $o\in M$ and write $V(R)=\operatorname{vol}(B(o,R))$ for geodesic ball volume. Let $0<s<1$ and $q>1$. If
\begin{equation}\label{eq:volume-manifold}
 \int_1^\infty\frac{R^{2sq-1}}{V(R)^{q-1}}\,dR=\infty,
\end{equation}
then every nonnegative weak supersolution of $A^su\ge u^q$ is zero almost everywhere. Every finite Green integral supersolution is also zero almost everywhere.
\end{maintheorem}
\end{samepage}

On weighted graphs, the local theory leads to a corresponding volume series. For the ordinary normalized Laplacian, Gu, Huang, and Sun \cite{GuHuangSun2023} proved a critical power-logarithmic bound under a uniform positive lower bound on adjacent transition probabilities, using nonlinear test functions and estimates for ratios of neighboring solution values. Monticelli, Punzo, and Somaglia \cite{MonticelliPunzoSomaglia2025} treated inequalities with a potential by combining weighted volume growth with an upper bound on the Laplacian of a chosen distance. Gu, Hao, Huang, and Sun \cite[Corollary 1.7]{GuHaoHuangSun2026} established the full divergence criterion with exponent $2q-1$ through Green-current flow decomposition, a pathwise Hardy estimate, and relative capacity bounds.

Fractional graph operators can be described by semigroup and extension formulas. Ciaurri, Roncal, Stinga, Torrea, and Varona \cite{CiaurriEtAl2018} developed these representations on the one-dimensional lattice; Zhang, Lin, and Yang studied finite graphs \cite{ZhangLinYang2026} and stochastically complete locally finite graphs \cite{ZhangLinYang2025}. Delmotte's heat-kernel theory \cite{Delmotte1999} gives geometric estimates for reversible Markov chains under volume, Poincar\'e, and local ellipticity assumptions. A fractional power also introduces jumps between nonadjacent vertices. To express a criterion using the geometry of the original graph, we keep its edges in the horizontal extension energy and recover the fractional scale from the extension height. This gives the following result for the normalized Laplacian.

\begin{samepage}
\begin{maintheorem}\label{thm:intro-graph}\label{thm:graph}
Let $X$ be an infinite connected locally finite loopless undirected graph with finite positive symmetric edge weights $c_{xz}$. Set
\[
 \mu(x)=\sum_zc_{xz},\qquad
 Pv(x)=\frac1{\mu(x)}\sum_zc_{xz}v(z),\qquad A=I-P.
\]
Fix $o\in X$ and put $V(R)=\sum_{d(o,x)\le R}\mu(x)$, where $d$ is the original combinatorial distance. Let $0<s<1$ and $q>1$. If
\begin{equation}\label{eq:volume}
 \sum_{n=1}^\infty\frac{n^{2sq-1}}{V(n)^{q-1}}=\infty,
\end{equation}
then every nonnegative pointwise or weak supersolution of $A^su\ge u^q$ is identically zero. Every finite Green integral supersolution is also identically zero.
\end{maintheorem}
\end{samepage}

Theorem~\ref{thm:graph} gives the fractional exponent $2sq-1$ with the normalization and original distance of the local graph criterion in \cite{GuHaoHuangSun2026}. Both of our volume conditions have the same consequences for polynomial growth. If $V(R)\le CR^D$ with $D>2s$, they exclude nonzero supersolutions for $1<q\le D/(D-2s)$; if $D\le2s$, they do so for every $q>1$. Corollary~\ref{cor:growth} records the critical logarithmic bound as well. The integral and series formulations also apply to irregular volume growth.

The extension argument is related to the author's recurrence criterion \cite{Cheng2026Recurrence}. That work proves recurrence of $A^s$ on manifolds when $\int_1^\infty R^{2s-1}/V(R)\,dR$ diverges, using a weighted trace estimate and radial cutoffs. Here the same extension weight is used to estimate the $L^q$ mass of a fixed-source Green potential. The nonlinear estimate introduces the exponent $q$ and leads to the volume conditions in Theorems~\ref{thm:manifold} and~\ref{thm:graph}.

The proof compares two estimates for $g=A_\Omega^{-s}f$, with a fixed bounded nonnegative source $f$ of positive mass and an expanding family of domains $\Omega$. We kill the original Laplacian in $\Omega$ and then form its fractional Green operator; Kim, Song, and Vondra\v{c}ek \cite{KimSongVondracek2019} discuss this order of killing and subordination. Nonlinear testing gives an $L^q$ upper bound independent of $\Omega$. The positive-source Stinga--Torrea extension gives a lower bound through horizontal annular energy and the decrease between successive heights. Bounded overlap and telescoping combine the scales, and the volume divergence contradicts the uniform upper bound. Gradients on manifolds and edge differences on graphs supply the horizontal estimates; the height-variable argument is shared.

\section{Operator conventions and the nonlinear Green estimate}\label{sec:upper}
We specify the operators and supersolutions used in Theorems~\ref{thm:manifold} and~\ref{thm:graph}. Throughout this section, $0<s<1$ and $q>1$.

\label{sec:manifold}
Let $M$ be a smooth complete connected Riemannian manifold without boundary, of any positive dimension. Let $\mu$ be Riemannian volume and let $A=-\Delta$ be the nonnegative Friedrichs realization on $L^2(M,\mu)$. Its heat semigroup $T_t=e^{-tA}$ is the minimal, positivity-preserving, sub-Markov semigroup. We do not assume stochastic completeness. For $0<s<1$ and $\psi\in C_c^\infty(M)$, its spectral fractional power is represented by
\begin{equation}\label{eq:semigroup-power}
 A^s\psi=c_s\int_0^\infty(\psi-T_t\psi)\,\frac{dt}{t^{1+s}},
 \qquad c_s=\frac{s}{\Gamma(1-s)}>0.
\end{equation}
If $\psi\ge0$, this function is nonpositive outside $\operatorname{supp}\psi$ and bounded on $M$. Indeed,
\[
 \norm{\psi-T_t\psi}_\infty\le t\norm{A\psi}_\infty\quad(0<t\le1),
 \qquad \norm{\psi-T_t\psi}_\infty\le2\norm\psi_\infty\quad(t\ge1).
\]
These estimates also fix the meaning of the nonlocal pairing below.

\begin{definition}\label{def:weak}
A manifold weak supersolution is a nonnegative $u\in L^1_{\mathrm{loc}}(M)$ such that, for every nonnegative $\psi\in C_c^\infty(M)$,
\begin{equation}\label{eq:weak-super}
 \int_M uA^s\psi\,d\mu\ge\int_Mu^q\psi\,d\mu.
\end{equation}
The left side is interpreted as an extended integral whose positive part is finite. The right side is initially allowed to be infinite.
\end{definition}
The inequality forces both $\int u|A^s\psi|<\infty$ and $\int u^q\psi<\infty$. Choosing $\psi\ge1$ on an arbitrary compact set gives $u\in L^q_{\mathrm{loc}}(M)$. Thus this definition does not impose global $L^q$ integrability or an additional tail hypothesis.

Fix $o\in M$. With geodesic distance $d$, write
\[
 B_R=\{x\in M:d(o,x)\le R\},\qquad V(R)=\mu(B_R).
\]
These closed balls are compact. For a nonnegative measurable function $v$ on $M$, define
\begin{equation}\label{eq:global-green}
 \mathcal G_sv(x)=\frac1{\Gamma(s)}\int_0^\infty t^{s-1}T_tv(x)\,dt\in[0,\infty].
\end{equation}
A manifold Green integral supersolution is a nonnegative measurable function $u$, finite almost everywhere, such that $u\ge\mathcal G_s(u^q)$ almost everywhere. This definition allows the Green operator to have extended values and does not presuppose a bounded inverse or transience.

We first dispose of compact manifolds and a local integrability issue for integral supersolutions. The minimal heat kernel on a connected manifold is continuous and strictly positive for positive time. If $u\ge\mathcal G_s(u^q)$ is finite almost everywhere, choose a point $x_0$ where the inequality holds with $u(x_0)<\infty$. For each compact $K$, the heat kernel has a positive minimum on $\{x_0\}\times K\times[1,2]$. Formula \eqref{eq:global-green} therefore implies $\int_Ku^q<\infty$. On a compact manifold, weak supersolutions vanish by testing $\psi=1$. Integral supersolutions belong to $L^q(M)\subset L^1(M)$; conservation of the compact heat semigroup and Tonelli give
\[
 \int_Mu\,d\mu\ge\frac1{\Gamma(s)}\int_0^\infty t^{s-1}\,dt\int_Mu^q\,d\mu,
\]
which forces $u=0$. Hence all subsequent manifold domains lie in a noncompact $M$.

\label{sec:graph}
For the graph case, let $X$ be an infinite connected locally finite undirected graph without loops, and let $E$ be its set of unordered edges. Assign finite positive symmetric weights $c_{xz}=c_{zx}$ to its edges and set $c_{xz}=0$ otherwise. We also write $c_e=c_{xz}$ for $e=\{x,z\}\in E$. Define
\begin{equation}\label{eq:setup}
 \mu(x)=\sum_zc_{xz},\qquad
 Pv(x)=\frac1{\mu(x)}\sum_zc_{xz}v(z),\qquad A=I-P.
\end{equation}
Then $P$ is a self-adjoint Markov contraction on $\ell^2(X,\mu)$, and the spectrum of $A$ lies in $[0,2]$. The positive kernel representation
\begin{equation}\label{eq:subordination}
 A^s=I-K_s,\qquad K_s=\sum_{j\ge1}a_jP^j,
 \quad a_j=(-1)^{j+1}\binom{s}{j}>0,\quad \sum_{j\ge1}a_j=1
\end{equation}
follows from the binomial series. Absolute summability extends the scalar identity to both endpoints of the spectrum of $P$ and gives operator-norm convergence on $\ell^2$.

A graph pointwise supersolution is a finite function $u:X\to[0,\infty)$ with $K_su(x)<\infty$ and
\begin{equation}\label{eq:supersolution}
 u(x)-K_su(x)\ge u(x)^q\qquad(x\in X).
\end{equation}
It automatically satisfies $0\le u\le1$. The weak formulation using nonnegative finitely supported tests and absolutely convergent pairings is equivalent: testing at one vertex and using reversibility yields \eqref{eq:supersolution}; conversely boundedness of $u$ and boundedness of $A^s=I-K_s$ on $\ell^1$ justify the pairings. Graph heat operators are given by $T_t=e^{-t}\sum_{j\ge0}t^jP^j/j!$.

Fix $o\in X$ and let $d$ be the original combinatorial distance, the minimum number of edges in a path. Put
\[
 B_R=\{x\in X:d(o,x)\le R\},\qquad V(R)=\sum_{x\in B_R}\mu(x).
\]
The balls are finite by local finiteness. For $v:X\to[0,\infty]$, define
\begin{equation}\label{eq:global-green-graph}
 \mathcal G_sv(x)=\frac1{\Gamma(s)}\int_0^\infty t^{s-1}T_tv(x)\,dt\in[0,\infty].
\end{equation}
A graph Green integral supersolution is a finite nonnegative function satisfying $u(x)\ge\mathcal G_s(u^q)(x)$ at every vertex. All volumes and operators in the graph case use the original measure $\mu(x)=\sum_zc_{xz}$.

\label{sec:common}
For the estimates below, fix either the manifold setting of Theorem~\ref{thm:manifold} or the graph setting of Theorem~\ref{thm:graph}. Write $X=M$ in the manifold case and retain $X$ for the vertex set in the graph case. The measure $\mu$, distance $d$, balls $B_R$, volume $V(R)$, operator $A$, and Green operator $\mathcal G_s$ always belong to that fixed setting. Integrals against $\mu$ are ordinary volume integrals on a manifold and weighted sums on a graph. Put
\[
 p=\frac q{q-1},\qquad 0<s<1,\quad q>1.
\]
For polynomial volume growth, the criterion gives an explicit range of exponents. The logarithmic factor below also records the borderline case.

\begin{samepage}
\begin{corollary}\label{cor:growth}
In either setting assume, for all sufficiently large $R$,
\[
 V(R)\le C R^D(\log R)^\gamma,
\]
where $D\ge0$ and $\gamma\in\mathbb R$. The nonexistence conclusion of Theorems~\ref{thm:manifold} and~\ref{thm:graph} holds if
\[
 D(q-1)<2sq,
 \qquad\text{or if}\qquad D(q-1)=2sq\ \text{and}\ \gamma(q-1)\le1.
\]
In particular, a polynomial bound $V(R)\le CR^D$ excludes nonzero supersolutions for $1<q\le D/(D-2s)$ when $D>2s$, and for all $q>1$ when $D\le2s$.
\end{corollary}
\end{samepage}
\begin{proof}
The integrand or summand in the volume condition is bounded below, for large radii, by a positive constant times
\[
 R^{2sq-1-D(q-1)}(\log R)^{-\gamma(q-1)}.
\]
A power greater than $-1$ gives divergence for every fixed logarithmic power. At exponent $-1$, the integral and series diverge exactly when $\gamma(q-1)\le1$. For a pure polynomial bound, $D(q-1)\le2sq$ rearranges as stated when $D>2s$, and is strict for all $q>1$ when $D\le2s$.
\end{proof}

Our next aim is an upper bound that is uniform over the finite domain. The source will be chosen once from a putative supersolution and then kept fixed as the domain grows. We first specify the killed operators, since the order of killing and taking a fractional power determines the comparison used in the nonlinear test.

On a manifold, $\Omega\Subset M$ is a smooth connected domain and
$A_\Omega=-\Delta_\Omega^D$ is the ordinary Dirichlet Laplacian. The domainwise Dirichlet Poincar\'e inequality gives a positive first eigenvalue, so $A_\Omega$ is strictly positive and self-adjoint. For a nonnegative test $\Psi$ compactly supported in $\Omega$, killed heat domination and \eqref{eq:semigroup-power} give
\begin{equation}\label{eq:test-comparison}
 A^s\widetilde\Psi\le A_\Omega^s\Psi\ \text{on }\Omega,
 \qquad A^s\widetilde\Psi\le0\ \text{off }\Omega.
\end{equation}
The tilde denotes zero extension to $M$.

On a graph, $\Omega$ is a finite connected vertex set, and
\begin{equation}\label{eq:killed}
 P_\Omega v(x)=\sum_{z\in\Omega}\frac{c_{xz}}{\mu(x)}v(z),\qquad A_\Omega=I-P_\Omega.
\end{equation}
The denominator is the original $\mu(x)$; exit is killed without renormalization. For a function $v$ and its zero extension $\widetilde v$ to $X$,
\begin{equation}\label{eq:energy}
 \ip{A_\Omega v}{v}=\sum_{\{x,z\}\in E}c_{xz}(\widetilde v(x)-\widetilde v(z))^2.
\end{equation}
Connectedness and a nonempty exterior make this positive unless $v=0$, so $A_\Omega$ is strictly positive and self-adjoint. Paths staying in $\Omega$ form a subset of all paths. Hence for nonnegative $u$,
\begin{equation}\label{eq:comparison}
 P_\Omega^j(u|_\Omega)\le(P^ju)|_\Omega,
 \qquad A_\Omega^s(u|_\Omega)\ge(A^su)|_\Omega.
\end{equation}

An admissible domain means the domain just specified in the fixed setting. Write
\[
 T_t^\Omega=e^{-tA_\Omega},\qquad G_{\Omega,s}=A_\Omega^{-s}.
\]
The killed heat semigroup is positive and sub-Markov, with a spectral gap, and
\begin{equation}\label{eq:finite-green}
 G_{\Omega,s}h=\frac1{\Gamma(s)}\int_0^\infty t^{s-1}T_t^\Omega h\,dt.
\end{equation}
The Green operator preserves nonnegativity and improves positivity on the connected domain. For bounded $h$, its potential is bounded: small times use $L^\infty$ contraction; at large times the spectral gap and a fixed-time $L^2$ to $L^\infty$ heat bound give exponential decay. This fixed-time bound follows from the bounded positive-time Dirichlet heat kernel on a smooth relatively compact manifold domain, and from finite dimensionality on a graph. Constants used only for this boundedness may depend on $\Omega$.

We use the spectral Green operator, formed by killing before taking a fractional power. It need not equal the inverse of the part operator $v\mapsto(A^s\widetilde v)|_\Omega$. The geometry-specific comparisons are \eqref{eq:test-comparison} on manifolds and \eqref{eq:comparison} on graphs. In either setting, positivity and killed heat domination also give
\begin{equation}\label{eq:green-comparison}
 G_{\Omega,s}(v|_\Omega)\le(\mathcal G_sv)|_\Omega\qquad(v\ge0).
\end{equation}

\begin{lemma}\label{lem:upper}
Let $u$ be a nonzero supersolution in the fixed setting, as in Theorem~\ref{thm:manifold} or Theorem~\ref{thm:graph}. Choose a measurable set $F$ contained in a compact set, finite in the graph case, such that
\[
 f=u^q\one_F\quad\text{is bounded},\qquad m=\int_X f\,d\mu>0.
\]
Then every admissible $\Omega$ containing that compact set satisfies
\begin{equation}\label{eq:green-upper}
 \norm{G_{\Omega,s}f}_{L^q(\Omega,\mu)}^q\le p^q m.
\end{equation}
\end{lemma}
\begin{proof}
Such an $F$ exists by local $L^q$ integrability and nontriviality. For a graph every bounded finite set of values is sufficient. Set $L=A_\Omega^s$. In the manifold case take $0\le h\in C_c^\infty(\Omega)$; in the graph case take arbitrary $h\ge0$ on $\Omega$. Write $\phi=G_{\Omega,s}h$ and set
\[
 T=\int_\Omega u^q\phi^p\,d\mu.
\]
We will first prove that $T$ is finite and satisfies \eqref{eq:T-bound}, then use H\"older and duality to obtain \eqref{eq:green-upper}.

The convex chain inequality for a smooth convex $F_0$ with $F_0(0)=0$ is
\begin{equation}\label{eq:convex-smooth}
 L F_0(\phi)\le F_0'(\phi)L\phi.
\end{equation}
Indeed, apply Jensen to $T_t^\Omega$ after placing its missing mass at value zero, and then the tangent inequality at $\phi$; integrate $F_0(\phi)-T_t^\Omega F_0(\phi)\le F_0'(\phi)(\phi-T_t^\Omega\phi)$ against $c_st^{-1-s}dt$. For the functions below all terms are defined. Our target for each solution class is
\begin{equation}\label{eq:T-bound}
 T\le p\int_\Omega u\phi^{p-1}h\,d\mu<\infty.
\end{equation}

\smallskip\noindent\textit{Manifold weak supersolutions.}
A smooth cutoff of the power is needed. The function $\phi$ is smooth on $\overline\Omega$ and vanishes on $\partial\Omega$: $h\in D(A_\Omega^k)$ for every integer $k$, so the spectral gap implies the same for $A_\Omega^{-s}h$, and elliptic boundary regularity applies. Choose a smooth nondecreasing $\chi:\mathbb R\to[0,1]$, zero on $(-\infty,1]$ and one on $[2,\infty)$, and put
\begin{equation}\label{eq:F-eps}
 F_\varepsilon(t)=p\int_0^t r^{p-1}\chi(r/\varepsilon)\,dr\quad(t\ge0),
 \qquad F_\varepsilon(t)=0\quad(t<0).
\end{equation}
It is smooth and convex, is zero near zero, and satisfies
\[
 F_\varepsilon(t)\uparrow t^p,\qquad
 0\le F_\varepsilon'(t)\le pt^{p-1},\qquad
 F_\varepsilon'(t)\longrightarrow pt^{p-1}\quad(t\ge0).
\]
Thus $F_\varepsilon(\phi)\in C_c^\infty(\Omega)$ and its zero extension is a legitimate global test for every $p>1$. Using \eqref{eq:test-comparison} and \eqref{eq:convex-smooth} gives
\[
 \int_\Omega u^qF_\varepsilon(\phi)\,d\mu
 \le\int_\Omega uF_\varepsilon'(\phi)h\,d\mu.
\]
The right side is dominated by $p u\phi^{p-1}h$, which is integrable because $h$ has compact support and $\phi$ is bounded. Monotone convergence on the left and dominated convergence on the right give \eqref{eq:T-bound}, and in particular $T<\infty$, before any use of H\"older.

\smallskip\noindent\textit{Graph pointwise supersolutions.}
On a graph the same argument with the substochastic matrix $Q=\sum_ja_jP_\Omega^j$ also gives directly
\begin{equation}\label{eq:convex}
 L(\phi^p)\le p\phi^{p-1}L\phi.
\end{equation}

For a graph pointwise supersolution, \eqref{eq:comparison}, self-adjointness and \eqref{eq:convex} yield
\[
 T\le\ip{Lu}{\phi^p}=\ip u{L(\phi^p)}
 \le p\int_\Omega u\phi^{p-1}h\,d\mu.
\]
All sums here are finite.

\smallskip\noindent\textit{Green integral supersolutions in either setting.}
For an integral supersolution, \eqref{eq:green-comparison} gives $G_{\Omega,s}(u^q|_\Omega)\le u$. On a manifold apply $G_{\Omega,s}$ to \eqref{eq:convex-smooth} with $F_\varepsilon$ and then use positivity and monotone convergence. On a graph apply it to \eqref{eq:convex}. In both cases,
\[
 \phi^p\le pG_{\Omega,s}(h\phi^{p-1}).
\]
Tonelli and symmetry therefore give
\[
 T\le p\int_\Omega G_{\Omega,s}(u^q)h\phi^{p-1}\,d\mu
 \le p\int_\Omega uh\phi^{p-1}\,d\mu<\infty,
\]
where local integrability and the support of $h$ justify the last finiteness on a manifold.

\smallskip\noindent\textit{Common conclusion.}
In each case, \eqref{eq:T-bound} is now justified with $T<\infty$. Since $(p-1)q=p$, H\"older gives
\[
 T\le pT^{1/q}\norm h_p,
 \qquad \norm{G_{\Omega,s}h}_{L^p(\Omega,u^q\mu)}=T^{1/p}\le p\norm h_p.
\]
The conclusion is immediate if $T=0$; otherwise divide by $T^{1/q}$. Finally,
\begin{align}\label{eq:weighted-green}
 \ip{G_{\Omega,s}f}{h}
 &=\int_Fu^qG_{\Omega,s}h\,d\mu\notag\\
 &\le m^{1/q}\norm{G_{\Omega,s}h}_{L^p(u^q\mu)}
 \le p m^{1/q}\norm h_p.
\end{align}
Duality gives \eqref{eq:green-upper}. For a manifold, nonnegative compactly supported smooth functions suffice by density in $L^p(\Omega)$.
\end{proof}

\section{A multiscale lower bound}\label{sec:lower}
The upper bound in Lemma~\ref{lem:upper} depends on the supersolution. The geometric estimate below applies to every bounded nonnegative source of positive mass. It measures the $L^q$ cost of its Green potential over a finite collection of scales.

\begin{samepage}
\begin{proposition}\label{prop:lower}
In the fixed setting, with $M$ noncompact in the manifold case, let $f\ge0$ be bounded, nonzero, and compactly supported, and put $m=\int_Xf\,d\mu$. Choose an integer $R_0\ge16$ with $\operatorname{supp}f\subset B_{R_0/16}$ and set $R_k=2^kR_0$. For every $N\ge0$ and every admissible domain $\Omega$ containing $B_{16R_N}$ as above,
\begin{equation}\label{eq:green-lower}
 \norm{G_{\Omega,s}f}_{L^q(\Omega,\mu)}^q
 \ge c_{s,q}m^q\sum_{k=0}^N\frac{R_k^{2sq}}{V(8R_k)^{q-1}}.
\end{equation}
The positive constant depends only on $s,q$, not on the domain, dimension, graph, weights, source, or number of scales.
\end{proposition}
\end{samepage}

The proof will use one domain $\Omega$ for all the scales under consideration. We first extend $G_{\Omega,s}f$ in a height variable. An energy estimate controls the horizontal flux through spatial annuli, and the decrease in height controls the vertical flux. Their sum bounds the contribution of each scale in \eqref{eq:green-lower}.

\subsection{The positive-source extension}\label{sec:extension}
For an admissible $\Omega$ and bounded nonzero nonnegative $f$ supported in a compact subset of $\Omega$, let $g=G_{\Omega,s}f$ and define
\begin{equation}\label{eq:extension}
 U(x,y)=\frac1{\Gamma(s)}\int_0^\infty t^{s-1}e^{-y^2/(4t)}T_t^\Omega f(x)\,dt,
 \qquad y\ge0.
\end{equation}
Put $a=1-2s$ and
\begin{equation}\label{eq:ds}
 d_s=2^{1-2s}\frac{\Gamma(1-s)}{\Gamma(s)}.
\end{equation}
The horizontal Dirichlet form, with zero extension understood on a graph, is
\begin{equation}\label{eq:horizontal-form}
 \cE_\Omega(v,w)=
 \begin{cases}
 \displaystyle\int_\Omega\langle\nabla v,\nabla w\rangle\,d\mu,&\text{on a manifold},\\[4pt]
 \displaystyle\sum_{\{x,z\}\in E}c_{xz}D_evD_ew,&\text{on a graph}.
 \end{cases}
\end{equation}
Here $D_ev=v(x)-v(z)$ for one arbitrary orientation of each unordered edge. Both differences change sign together if the orientation is reversed.

\begin{lemma}\label{lem:extension}
The extension satisfies, in $L^2(\Omega)$,
\begin{gather}\label{eq:extension-pde}
 \partial_y(y^aU_y)=y^aA_\Omega U,\qquad U(0)=g,\\
 \label{eq:extension-flux}
 -\lim_{y\downarrow0}y^aU_y=d_sf.
\end{gather}
Moreover $U_y\le0$ and $0\le U(y_2)\le U(y_1)\le g$ for $y_2>y_1\ge0$, almost everywhere. It has finite weighted energy and obeys
\begin{equation}\label{eq:extension-energy}
 \int_0^\infty y^a\bigl(\norm{U_y}_2^2+\cE_\Omega(U,U)\bigr)\,dy
 =d_s\ip fg<\infty.
\end{equation}
For bounded finite-energy tests $\Psi$ with compact horizontal support in $\Omega$, compact support in $y\ge0$, and their energy trace at zero,
\begin{equation}\label{eq:weak-extension}
 \int_0^\infty y^a\left[\ip{U_y}{\Psi_y}+\cE_\Omega(U,\Psi)\right]dy
 =d_s\int_\Omega f\Psi(0)\,d\mu.
\end{equation}
On a graph the trace and flux limits are pointwise; on a manifold the asserted boundary limits are $L^2$ and weak limits, rather than everywhere classical limits for measurable $f$.
\end{lemma}
\begin{proof}
Since $f\in L^2$ and $A_\Omega$ has a positive spectral gap, $g\in D(A_\Omega^s)$. Formula \eqref{eq:extension} is the representation in \cite[(1.7)--(1.8)]{StingaTorrea2010} with boundary value $g$ and source $A_\Omega^sg=f$. Its basic identities can also be checked spectrally as follows. For $y>0$,
\[
 (\partial_{yy}+a y^{-1}\partial_y)e^{-y^2/(4t)}
 =\left(\frac{y^2}{4t^2}-\frac{1-s}{t}\right)e^{-y^2/(4t)}.
\]
After multiplication by $t^{s-1}$ this is the derivative in $t$ of $t^{s-1}e^{-y^2/(4t)}$. Integration by parts in $t$ gives \eqref{eq:extension-pde}; the factor at zero and the spectral gap at infinity remove the boundary terms. Positivity of the integrand gives the asserted monotonicity and $U\le g$, as well as the trace by convergence of the Green integral.

Differentiating and substituting $t=y^2/(4\tau)$ gives
\[
 -y^aU_y=\frac{2^{1-2s}}{\Gamma(s)}
 \int_0^\infty\tau^{-s}e^{-\tau}T_{y^2/(4\tau)}^\Omega f\,d\tau
 \longrightarrow d_sf\quad\text{in }L^2.
\]
The bound $\norm{T_t^\Omega f}_2\le\norm f_2$ and integrability of $\tau^{-s}e^{-\tau}$ justify this limit. The same calculation works pointwise on a finite graph.

For completeness, expand in Dirichlet eigenfunctions. A scalar coefficient with boundary value $g_j$ and eigenvalue $\lambda_j>0$ satisfies $(y^av_j')'=\lambda_j y^av_j$, decays at infinity, and has bottom flux $d_s\lambda_j^s g_j$. Multiplication by $v_j$ and integration give
\[
 \int_0^\infty y^a\bigl(|v_j'|^2+\lambda_j|v_j|^2\bigr)dy
 =d_s\lambda_j^s|g_j|^2.
\]
Summing nonnegative energies proves \eqref{eq:extension-energy}, because $\sum_j\lambda_j^s|g_j|^2=\ip fg$. The scalar decay and endpoint limits follow directly from the displayed heat integral; one may first integrate over $[\varepsilon,L]$.

Finite spectral sums give the weak identity, and energy convergence extends it to finite-energy tests by Cauchy--Schwarz and the trace inequality. In particular, all compactly supported Lipschitz cutoffs below and the bounded Lipschitz compositions used with $U$ are admissible. In the manifold setting this can equivalently be obtained by local weighted Sobolev approximation; $a\in(-1,1)$ ensures the usual energy traces. The boundary sign is positive since the contribution from $y=\varepsilon$ is $-\varepsilon^aU_y(\varepsilon)\Psi(\varepsilon)$.
\end{proof}

The sign $U_y\le0$ follows from the nonnegative source in \eqref{eq:extension} and will control the vertical flux. On a graph the horizontal coordinate remains discrete; the continuous variable $y$ is added at each vertex.

\subsection{Horizontal annular energy}\label{sec:energy}
We next estimate the energy through which the extension transports flux across a spatial annulus. The estimate must be uniform in the domain and must work for every $q>1$. Choosing the auxiliary exponent $r=\min(q,2)$ allows the manifold chain rule and a two-point graph inequality to produce the same bound.

Fix $R\ge16$, assume $\operatorname{supp}f\subset B_{R/16}$ and $B_{16R}\subset\Omega$, and set $r=\min(q,2)\in(1,2]$. In the manifold case the domain contains the closed ball as a compact subset. We use the same radial cutoffs in both settings. Let $\eta_R$ be one on $B_R$, linear in distance on $[R,2R]$, and zero beyond $B_{2R}$. Let $\theta_R$ be zero on $[0,R/4]$, rise linearly to one at $R/2$, equal one on $[R/2,3R]$, and decrease linearly to zero at $4R$, vanishing thereafter. Define
\begin{equation}\label{eq:annulus}
 \mathcal A_R=B_{8R}\setminus B_{R/8},\qquad
 \mathscr S_R=\{e:D_e\eta_R\ne0\}\quad\text{on a graph}.
\end{equation}
\smallskip\noindent\textit{Manifold cutoff bounds.}
We have $|\nabla\eta_R|\le R^{-1}$ and $|\nabla\theta_R|\le4R^{-1}$ almost everywhere.

\smallskip\noindent\textit{Graph cutoff bounds.}
The corresponding edge differences have the same bounds. Both endpoints of each edge in $\mathscr S_R$ lie where $\theta_R=1$. Every edge on which $\theta_R$ or its difference is active has both endpoints in $\mathcal A_R$: adjacent distances differ by at most one, and $R\ge16$ leaves the required slack. In both settings $\theta_R$ vanishes on the source.

Let $\zeta_R$ be one on $[0,2R]$, decrease linearly to zero on $[2R,4R]$, and vanish thereafter. Thus $|\zeta_R'|\le(2R)^{-1}$. The manifold energy calculation uses the ordinary chain rule. The next lemma is used only in the graph calculation; it supplies the edge inequality invoked below.

\begin{lemma}\label{lem:edge}
For $1<r\le2$, $v,w>0$ and $\theta,\eta\in[0,1]$,
\begin{align}\label{eq:edge}
 &(v-w)(v^{r-1}\theta^2-w^{r-1}\eta^2)\notag\\
 &\quad\ge\frac{r-1}{4}(v-w)^2(v+w)^{r-2}(\theta^2+\eta^2)
 -\frac2{r-1}(v^r+w^r)(\theta-\eta)^2.
\end{align}
\end{lemma}
\begin{proof}
Writing $S=\theta^2+\eta^2$, split the left side into
\[
 \frac12(v-w)(v^{r-1}-w^{r-1})S
 +\frac12(v-w)(v^{r-1}+w^{r-1})(\theta^2-\eta^2).
\]
The integral formula for the difference of powers, with $r-2\le0$, bounds the first term below by $(r-1)(v-w)^2(v+w)^{r-2}S/2$. The second has absolute value at most $\alpha\beta$, where
\[
 \alpha=|v-w|(v+w)^{(r-2)/2}S^{1/2},\qquad
 \beta=\frac{v^{r-1}+w^{r-1}}{\sqrt2}(v+w)^{(2-r)/2}|\theta-\eta|.
\]
Use $\alpha\beta\le(r-1)\alpha^2/4+\beta^2/(r-1)$. Concavity of $t^{r-1}$ and convexity of $t^r$ give
\[
 (v^{r-1}+w^{r-1})^2(v+w)^{2-r}
 \le2^{4-2r}(v+w)^r\le2^{3-r}(v^r+w^r).
\]
Since $2^{2-r}\le2$, the remaining error is bounded as in \eqref{eq:edge}.
\end{proof}

The quantity needed for horizontal flux is
\begin{equation}\label{eq:localized-energy}
 E_{r,R}=\begin{cases}
 \displaystyle\int_0^{2R}\!y^a\int_{B_{2R}\setminus B_R}
 U^{r-2}|\nabla_x U|^2\,d\mu\,dy,&\text{on a manifold},\\[5pt]
 \displaystyle\int_0^{2R}\!y^a\sum_{e\in\mathscr S_R}
 c_e(D_eU)^2(U_x+U_z)^{r-2}\,dy,&\text{on a graph}.
 \end{cases}
\end{equation}
All negative-power expressions are on positive values for $y>0$; limits at the bottom are understood through the energy integrals.

\begin{lemma}\label{lem:caccioppoli}
In either setting,
\begin{equation}\label{eq:caccioppoli}
 E_{r,R}\le C_{s,r}R^{-2s}\int_{\mathcal A_R}g^r\,d\mu.
\end{equation}
The constant depends only on $s,r$.
\end{lemma}
\begin{proof}
\smallskip\noindent\textit{Manifold calculation.}
Put $\Theta=\theta_R\zeta_R$ and test \eqref{eq:weak-extension} with $(U+\varepsilon)^{r-1}\Theta^2$, $\varepsilon>0$. This is an admissible energy test: $U$ is bounded, the derivative of the composition is bounded for fixed $\varepsilon$, and the cutoffs have compact support inside $\Omega$. Its source term is zero. With $\nabla$ denoting the full space--height gradient, expansion and Young's inequality yield
\begin{align*}
 &(r-1)\int y^a(U+\varepsilon)^{r-2}|\nabla U|^2\Theta^2
 =-2\int y^a(U+\varepsilon)^{r-1}\Theta\langle\nabla U,\nabla\Theta\rangle,\\
 &\int y^a(U+\varepsilon)^{r-2}|\nabla U|^2\Theta^2
 \le\frac4{(r-1)^2}\int y^a(U+\varepsilon)^r|\nabla\Theta|^2.
\end{align*}
For $r<2$ the left integrands increase as $\varepsilon\downarrow0$; the right side converges by domination with $(g+1)^r$ on a finite-measure cylinder. For $r=2$ the same estimate is the ordinary linear energy estimate. Restricting the left side to the horizontal energy in \eqref{eq:localized-energy} gives
\[
 E_{r,R}\le\frac{C_r}{R^2}\int_0^{4R}y^a\,dy\int_{\mathcal A_R}g^r\,d\mu.
\]
Only the Lipschitz gradient bound for distance is used.

\smallskip\noindent\textit{Graph calculation.}
Use $\Psi=U^{r-1}\theta_R^2\zeta_R^2$. For fixed $\Omega$ and $[0,4R]$, the positive continuous $U$ has a positive minimum on the finite vertex set, so the test is legitimate also for $r<2$. No uniform minimum is needed. Apply Lemma~\ref{lem:edge} to the horizontal form and the inequality
\[
 U_y\partial_y(U^{r-1}\zeta_R^2)
 \ge\frac{r-1}{2}U^{r-2}U_y^2\zeta_R^2
 -\frac2{r-1}U^r|\zeta_R'|^2
\]
to the vertical term. Retaining the desired horizontal energy and discarding the nonnegative vertical energy gives
\begin{align*}
 E_{r,R}\le C_r\int_0^{4R}y^a\bigg[
 &\zeta_R^2\sum_ec_e(U_x^r+U_z^r)(D_e\theta_R)^2\\
 &+|\zeta_R'|^2\sum_x\mu(x)U_x^r\theta_R(x)^2\bigg]dy.
\end{align*}
Edges touching the exterior have zero cutoff at both endpoints and contribute nothing. All other active endpoints belong to $\mathcal A_R$. Since $U\le g$ and $\sum_zc_{xz}=\mu(x)$,
\begin{equation}\label{eq:degree-absorb}
 \sum_ec_e(g_x^r+g_z^r)(D_e\theta_R)^2
 \le\frac{16}{R^2}\sum_{\mathcal A_R}\mu(x)g(x)^r.
\end{equation}
The vertical error has the same $R^{-2}$ bound.

\smallskip\noindent\textit{Common conclusion.}
In both settings,
\[
 R^{-2}\int_0^{4R}y^a\,dy=\frac{4^{2-2s}}{2-2s}R^{-2s},
\]
which proves the lemma with constants independent of $\Omega$ and the ambient geometry.
\end{proof}

\subsection{Flux and summation over scales}\label{sec:flux}
The annular energy estimate controls the horizontal part of the source flux. The remaining vertical part will be bounded by a decrease of $U$ over a height interval. These two contributions can be summed over scales for one fixed domain and one fixed potential.

Let $\xi_R$ equal one on $[0,R]$, decrease linearly to zero on $[R,2R]$, and vanish above $2R$. Test \eqref{eq:weak-extension} with $\Psi=\eta_R\xi_R$. Since $\eta_R=1$ on the source and $m=\int f\,d\mu$,
\begin{equation}\label{eq:flux}
 d_sm=H_R+Z_R,
\end{equation}
where
\begin{align*}
 H_R&=\int_0^{2R}y^a\xi_R\cE_\Omega(U,\eta_R)\,dy,\\
 Z_R&=\int_R^{2R}y^a\xi_R'\int_\Omega\eta_R U_y\,d\mu\,dy\ge0.
\end{align*}
The sign of $Z_R$ uses $U_y\le0$ and $\xi_R'\le0$; no sign is assigned to $H_R$. Define
\begin{equation}\label{eq:MN}
 \begin{split}
 W_R&=V(8R),\qquad M_R=\int_{\mathcal A_R}g^q\,d\mu,\\
 \delta_R(x)&=U(x,R)-U(x,2R)\ge0,\qquad
 N_R=\int_\Omega\delta_R^q\,d\mu.
 \end{split}
\end{equation}

\begin{lemma}\label{lem:single-scale}
For the same $g$ and $U$,
\begin{equation}\label{eq:single-scale}
 m^q\frac{R^{2sq}}{V(8R)^{q-1}}\le C_{s,q}(M_R+N_R).
\end{equation}
\end{lemma}
\begin{proof}
\smallskip\noindent\textit{Manifold horizontal flux.}
Weighted Cauchy--Schwarz gives
\[
 |H_R|\le E_{r,R}^{1/2}
 \left(\int_0^{2R}y^a\xi_R^2\int_{B_{2R}\setminus B_R}
 U^{2-r}|\nabla\eta_R|^2\,d\mu\,dy\right)^{1/2}.
\]
\smallskip\noindent\textit{Graph horizontal flux.}
The same Cauchy--Schwarz estimate has second factor
\[
 \left(\int_0^{2R}y^a\xi_R^2\sum_{e\in\mathscr S_R}
 c_e(D_e\eta_R)^2(U_x+U_z)^{2-r}\,dy\right)^{1/2}.
\]
For $r<2$, subadditivity of $t^{2-r}$ bounds $(U_x+U_z)^{2-r}$ by $g_x^{2-r}+g_z^{2-r}$; for $r=2$ the factor is one.

\smallskip\noindent\textit{Common flux bound.}
In either setting the squared second factor is at most
\[
 \frac{C}{R^2}\int_0^{2R}y^a\,dy
 \int_{\mathcal A_R}g^{2-r}\,d\mu
 \le C_sR^{-2s}\int_{\mathcal A_R}g^{2-r}\,d\mu.
\]
Here the graph estimate uses $\sum_zc_{xz}=\mu(x)$, as in \eqref{eq:degree-absorb}; the manifold estimate uses $|\nabla\eta_R|\le R^{-1}$. Combining this with Lemma~\ref{lem:caccioppoli} gives
\begin{equation}\label{eq:horizontal-product}
 |H_R|\le C_{s,r}R^{-2s}
 \left(\int_{\mathcal A_R}g^r\,d\mu\right)^{1/2}
 \left(\int_{\mathcal A_R}g^{2-r}\,d\mu\right)^{1/2}.
\end{equation}
At $r=2$ the second integral denotes $\mu(\mathcal A_R)$. For $0<t\le q$, H\"older gives
\[
 \int_{\mathcal A_R}g^t\,d\mu\le W_R^{1-t/q}M_R^{t/q}.
\]
Here $r\le q$ and $0\le2-r<q$. Thus the two exponents of $M_R$ sum to $r/(2q)+(2-r)/(2q)=1/q$, and the two volume exponents sum to $1-1/q$. Consequently
\begin{equation}\label{eq:horizontal}
 |H_R|\le C_{s,q}R^{-2s}W_R^{1-1/q}M_R^{1/q}.
\end{equation}

For the vertical flux, $y^a|\xi_R'|\le C_sR^{-2s}$ on $[R,2R]$. Since $U_y\le0$,
\begin{equation}\label{eq:vertical}
 Z_R\le C_sR^{-2s}\int_{B_{2R}}\delta_R\,d\mu
 \le C_sR^{-2s}W_R^{1-1/q}N_R^{1/q}.
\end{equation}
Combining $d_sm\le|H_R|+Z_R$ with \eqref{eq:horizontal}--\eqref{eq:vertical} and $(v+w)^q\le2^{q-1}(v^q+w^q)$ proves \eqref{eq:single-scale}.
\end{proof}

\begin{proof}[Proof of Proposition~\ref{prop:lower}]
Fix $N$ and one domain $\Omega$ large enough for all $k=0,\ldots,N$. All quantities below come from the same $g=G_{\Omega,s}f$ and its extension $U$. The annuli $\mathcal A_{R_k}$ have bounded overlap: for a positive distance $d(o,x)$, membership requires $d(o,x)/8\le R_k<8d(o,x)$, allowing at most seven dyadic scales. Therefore
\begin{equation}\label{eq:annulus-sum}
 \sum_{k=0}^NM_{R_k}\le7\norm g_q^q.
\end{equation}
The vertical increments are nonnegative and telescope pointwise almost everywhere:
\[
 \sum_{k=0}^N\delta_{R_k}(x)
 =U(x,R_0)-U(x,R_{N+1})\le g(x).
\]
Because $q>1$, $\sum_k\delta_{R_k}(x)^q\le(\sum_k\delta_{R_k}(x))^q\le g(x)^q$. Integration gives
\begin{equation}\label{eq:vertical-sum}
 \sum_{k=0}^NN_{R_k}\le\norm g_q^q.
\end{equation}
Summing Lemma~\ref{lem:single-scale} and applying \eqref{eq:annulus-sum}--\eqref{eq:vertical-sum} proves \eqref{eq:green-lower}.
\end{proof}

\section{Completion of the two proofs}\label{sec:sum}
We now combine the domain-independent upper bound with the multiscale lower bound. The only remaining issue is the factor $8$ in $V(8R)$: an exact shift of the dyadic index removes it without a doubling assumption.

\begin{proof}[Proof of Theorems~\ref{thm:manifold} and~\ref{thm:graph}]
The compact manifold case was proved in Section~\ref{sec:upper}. Otherwise fix the setting and suppose a nonzero supersolution exists, and fix $F,f,m$ as in Lemma~\ref{lem:upper}. Choose $R_0$ as in Proposition~\ref{prop:lower}. For each $N$, choose an admissible $\Omega$ containing $B_{16R_N}$. In the manifold proof, completeness makes the ball compact, so it lies in a smooth relatively compact connected domain. In the graph proof, local finiteness makes a sufficiently large connected ball finite. Combining the upper and lower estimates for the same $G_{\Omega,s}f$ gives
\begin{equation}\label{eq:contradiction-bound}
 m^{q-1}\sum_{k=0}^N\frac{R_k^{2sq}}{V(8R_k)^{q-1}}\le C_{s,q}p^q.
\end{equation}
The right side is independent of $N$ and $\Omega$.

Write $\alpha=2sq>0$, $\beta=q-1>0$.

\smallskip\noindent\textit{Manifold volume condition.}
For Theorem~\ref{thm:manifold}, volume monotonicity gives
\[
 \int_R^{2R}\frac{t^{\alpha-1}}{V(t)^\beta}\,dt
 \le C_\alpha\frac{R^\alpha}{V(R)^\beta}.
\]
\smallskip\noindent\textit{Graph volume condition.}
For Theorem~\ref{thm:graph}, the corresponding estimate for integer $R\ge1$ is
\[
 \sum_{n=R}^{2R-1}\frac{n^{\alpha-1}}{V(n)^\beta}
 \le C_\alpha\frac{R^\alpha}{V(R)^\beta}.
\]
\smallskip\noindent\textit{Common contradiction.}
In the chosen setting, its divergence assumption therefore implies $\sum_{k\ge0}R_k^\alpha/V(R_k)^\beta=\infty$. Since $8R_k=R_{k+3}$,
\begin{equation}\label{eq:dyadic-shift}
 \sum_{k=0}^\infty\frac{R_k^\alpha}{V(8R_k)^\beta}
 =8^{-\alpha}\sum_{k=3}^\infty\frac{R_k^\alpha}{V(R_k)^\beta}=\infty.
\end{equation}
This exact index shift requires no doubling condition. The unbounded partial sums contradict \eqref{eq:contradiction-bound} because $m>0$.
\end{proof}

\section*{Acknowledgments}
The author thanks Professor Bobo Hua for his support.

\section*{Use of AI tools}
The author used OpenAI's Codex to investigate candidate constructions,
check calculations, and solicit separate AI reviews of draft proofs. The
author determined the direction of the work, verified all mathematical
claims and arguments, prepared the final manuscript, and assumes full
responsibility for its contents.

\bibliographystyle{amsplain}
\bibliography{references}
\end{document}